\documentclass[12pt,english,serif,mathsf]{article}
\usepackage{mathptmx}

\usepackage[T1]{fontenc}
\usepackage[a4paper]{geometry}
\usepackage{fancyhdr}
\usepackage{color}
\usepackage{babel}
\usepackage{array}
\usepackage{longtable}
\usepackage{amsmath}
\usepackage{amsthm}
\usepackage{amssymb}
\usepackage{setspace}
\usepackage{esint}
\usepackage[unicode=true,
 bookmarks=false,
 breaklinks=false,pdfborder={0 0 0},backref=false,colorlinks=true]
 {hyperref}
\hypersetup{
 xetex,linkcolor=blue,citecolor=[rgb]{0,0.6,0},urlcolor=blue}

\makeatletter

\providecommand{\tabularnewline}{\\}

\theoremstyle{plain}
\newtheorem{thm}{\protect\theoremname}
  \theoremstyle{plain}
  \newtheorem{lem}{\protect\lemmaname}
  \theoremstyle{remark}
  \newtheorem{rem}{\protect\remarkname}

\usepackage{setspace} 
\usepackage{mathptmx} 

\makeatother

  \providecommand{\lemmaname}{Lemma}
  \providecommand{\remarkname}{Remark}
\providecommand{\theoremname}{Theorem}

\begin{document}
\noindent \textbf{\footnotesize{}Journal of Contemporary Mathematical
Analysis}{\footnotesize \par}

\noindent {\footnotesize{}(Izvestiya Natsionalnoi Akademii Nauk Armenii.
Matematika)}{\footnotesize \par}

\noindent \textbf{\footnotesize{}vol. 40, no. 2 (2005), pp. 14-27}{\footnotesize \par}

\title{Inverse Scattering Problem\\for Sturm--Liouville Operators}
\author{Hayk Asatryan\\Department of Mathematics\\Yerevan State University, Armenia}
\date{}
{\let\newpage\relax\maketitle}

\author{{\small{}}%
\begin{longtable}{p{14cm}}
\textbf{\small{}Abstract.}{\small{} On the space $L^{2}(\mathbb{R})$
the Sturm-Liouville operator $L$ with certain behavior of the potential
at infinity is considered. It is proved that $L$ is uniquely determined
by its scattering data. The recovery of $L$ is reduced to the solving
of a certain linear integral equation.}\tabularnewline
\end{longtable}}

\maketitle

\lhead{\noindent \emph{\footnotesize{}H. A. Asatryan}}

\rhead{\noindent \emph{\footnotesize{}Inverse Scattering Problem for Sturm--Liouville
Operators}}

\noindent \bigskip{}

\noindent Consider the differential operation (expression) $l(y)=-y''+qy$
on $\mathbb{R}=(-\infty,+\infty)$, where the coefficient (the potential)
$q$ is a real-valued measurable function of $x\in\mathbb{R}$, such
that 
\begin{equation}
\int\limits _{-\infty}^{0}|q(x)-a^{-}|dx+\int\limits _{0}^{\infty}|q(x)-a^{+}|dx<\infty
\end{equation}
with some constants $a^{\pm}\in\mathbb{R}$.

We say that $l$ is applicable to a function $y:\mathbb{R}\longmapsto\mathbb{C}$,
if $y$ has absolutely continuous first derivative in any closed interval
$[\alpha,\beta]\subset\mathbb{R}$. We set 
\[
[y(x),z(x)]=i(y(x)\overline{z'(x)}-y'(x)\overline{z(x)})
\]
for any functions $y$ and $z$ in $\mathbb{R}$, to which $l$ is
applicable. Besides, denote
\[
\mu_{0}=-\infty,\quad\mu_{1}=\min\{a^{+},a^{-}\},\quad\mu_{2}=\max\{a^{+},a^{-}\},\quad\mu_{3}=\infty,
\]
\[
\lambda_{j}^{\pm}(\mu)=(-1)^{j-1}\sqrt{\mu-a^{\pm}}\quad(\mu\in\mathbb{C},\ j=1,2)
\]
(where we take the principal value of the root). For $\mu\in\mathbb{R}$,
by $r^{\pm}(\mu)$ we denote the half of the number of real roots
of the equation $\lambda^{2}+a^{\pm}=\mu$. Obviously, the functions
$r^{\pm}\left(\cdot\right)$ are constant in each interval $(\mu_{k},\mu_{k+1})$,
$k=0,1,2$. For each $k=0,1,2$, by $r_{k}^{\pm}$ we denote the value
of the function $r^{\pm}$ in the interval $(\mu_{k},\mu_{k+1})$.
\begin{thm}
For any $\mu\in\mathbb{C}\setminus\{a^{+}\}$ (or $\mu\in\mathbb{C}\setminus\{a^{-}\}$)
the equation $l(y)=\mu y$ has linearly independent solutions $y_{1}^{+}(x,\mu)$,
$y_{2}^{+}(x,\mu)$ (correspondingly, $y_{1}^{-}(x,\mu)$, $y_{2}^{-}(x,\mu)$)
such that they and their derivatives with respect to $x$ have the
following asymptotics as $x\to\infty$ (correspondingly, as $x\to-\infty$):
\begin{equation}
y_{j}^{\pm}(x,\mu)=e^{ix\lambda_{j}^{\pm}(\mu)}[1+o(1)],
\end{equation}
\begin{equation}
y_{j}^{\prime\pm}(x,\mu)=i\lambda_{j}^{\pm}(\mu)e^{ix\lambda_{j}^{\pm}(\mu)}[1+o(1)].
\end{equation}
Moreover:
\begin{enumerate}
\item for any $\mu\in\mathbb{R}\setminus\{a^{+},a^{-}\}$ and $j,k=1,2$
\[
[y_{j}^{\pm}(x,\mu),y_{k}^{\pm}(x,\mu)]=\begin{cases}
2\lambda_{j}^{\pm}(\mu) & \mbox{\textup{if}}\quad\lambda_{j}^{\pm}(\mu)=\overline{\lambda_{k}^{\pm}(\mu)}\\
0 & \mbox{\textup{if}}\quad\lambda_{j}^{\pm}(\mu)\ne\overline{\lambda_{k}^{\pm}(\mu)}.
\end{cases}
\]

\item for each $k=0,1,2$, the functions $y_{j}^{+}(x,\mu)$, $y_{j}^{\prime+}(x,\mu)$
$(1\leqslant j\leqslant1+r_{k}^{+})$ and $y_{j}^{-}(x,\mu)$, $y_{j}^{\prime-}(x,\mu)$
$(2-r_{k}^{-}\leqslant j\leqslant2)$ are continuous in the variables
$x\in\mathbb{R}$, $\mu\in(\mu_{k},\mu_{k+1})$, 
\item the functions $y_{1}^{+}(x,\mu)$, $y_{1}^{\prime+}(x,\mu)$, $y_{2}^{-}(x,\mu)$,
$y_{2}^{\prime-}(x,\mu)$ are differentiable in $\mu\in(\mu_{0},\mu_{1})$,
and the derivatives $\frac{\partial y_{1}^{+}(x,\mu)}{\partial\mu},\frac{\partial y_{1}^{\prime+}(x,\mu)}{\partial\mu}$,
$\frac{\partial y_{2}^{-}(x,\mu)}{\partial\mu},\frac{\partial y_{2}^{\prime-}(x,\mu)}{\partial\mu}$
are continuous functions of two variables $x\in\mathbb{R}$, $\mu\in(\mu_{0},\,\mu_{1})$, 
\item if the potential $q$ satisfies the condition 
\begin{equation}
\int\limits _{-\infty}^{0}(1-x)|q(x)-a^{-}|\;dx+\int\limits _{0}^{\infty}(1+x)|q(x)-a^{+}|\;dx<\infty,
\end{equation}
then for each $k=1,2$ the functions $y_{j}^{+}(x,\mu)$, $y_{j}^{\prime+}(x,\mu)$
$(1\leqslant j\leqslant1+r_{k}^{+})$ and $y_{j}^{-}(x,\mu),\,y_{j}^{\prime-}(x,\mu)$
$(2-r_{k}^{-}\leqslant j\leqslant2)$ are differentiable in $\mu\in(\mu_{k},\mu_{k+1})$,
and the derivatives $\frac{\partial y_{j}^{+}(x,\mu)}{\partial\mu}$,
$\frac{\partial y_{j}^{\prime+}(x,\mu)}{\partial\mu}$ $(1\leqslant j\leqslant1+r_{k}^{+})$
as well as $\frac{\partial y_{j}^{-}(x,\mu)}{\partial\mu}$, $\frac{\partial y_{j}^{\prime-}(x,\mu)}{\partial\mu}$
$(2-r_{k}^{-}\leqslant j\leqslant2)$ are continuous in the variables
$x\in\mathbb{R}$, $\mu\in(\mu_{k},\mu_{k+1})$. 
\end{enumerate}
\end{thm}
\noindent \textit{Proof:} We outline it for $y_{j}^{+}(x,\mu)$ $(j=1,2)$.
For ${\hbox{\rm Im}}\lambda_{j}^{+}(\mu)\geqslant0$ the solutions
$y_{j}^{+}(x,\mu)$ are obtained from the integral equations 
\[
y_{j}^{+}(x,\mu)=e^{i\lambda_{j}^{+}(\mu)x}+\frac{1}{\lambda_{j}^{+}(\mu)}\int\limits _{x}^{\infty}\sin\left[\lambda_{j}^{+}(\mu)(t-x)\right]\,\left[q(t)-a^{+}\right]y_{j}^{+}(t,\mu)\;dt,
\]
and, if ${\hbox{\rm Im}}\lambda_{j}^{+}(\mu)<0$, from the integral
equations 
\[
y_{j}^{+}(x,\mu)=e^{i\lambda_{j}^{+}(\mu)x}-\frac{1}{2i\lambda_{j}^{+}(\mu)}\int\limits _{a}^{x}e^{-i\lambda_{j}^{+}(\mu)(x-t)}\left[q(t)-a^{+}\right]y_{j}^{+}(t,\mu)\;dt-
\]
\[
-\frac{1}{2i\lambda_{j}^{+}(\mu)}\int\limits _{x}^{\infty}e^{-i\lambda_{j}^{+}(\mu)(t-x)}\left[q(t)-a^{+}\right]y_{j}^{+}(t,\mu)\;dt,
\]
where the number $a$ is chosen to be large enough. These equations
are solved using the method of successive approximations. The relations
(2), (3) and the assertion 1 are obtained from these integral equations,
and the remaining assertions are proved using the first equation.
\qed

Note that for differential operators of order $m\geqslant2$ the existence
of solutions $y_{j}^{\pm}(x,\mu)$ $(j=1,2)$ satisfying the asymptotic
relations (2), (3) and the assertion 1) of Theorem 1, is established
in \cite{Pet03}.

The operator $L$ acting in the space $L^{2}(\mathbb{R})$ is defined
as follows (see \cite{Nai69}, p. 60). The domain $D$, where the
operator $L$ is defined, is the set of all functions $y\in L^{2}(\mathbb{R})$
for which the expression $l(y)$ is meaningful and $l(y)\in L^{2}(\mathbb{R})$,
and for any $y\in D$ we define $Ly=l(y)$. Under the condition (1),
the operator $L$ is self-adjoint (see \cite{Pet03}) and has a bounded
point spectrum. Moreover, under the condition (1), all eigenvalues
of $L$ (if they exist) are simple, lie in the interval $(-\infty,\mu_{1}]$
and as a limit point can have only $\mu_{1}$, see \cite{Asat06}.
In the same work, it is shown that under the condition (4) the set
of eigenvalues of the operator $L$ is finite and lies in the interval
$(-\infty,\mu_{1})$.

The inverse scattering problem for $L$ was considered by L. D. Faddeev
(see \cite{Fad64} and \cite{Mar77}, pp. 264 -- 283) for $a^{+}=a^{-}=0$,
and the case $a^{+}=a^{-}\ne0$ is easily the same. We shall assume
that $a^{+}\ne a^{-}$ and note that our approach is applicable to
the case $a^{+}=a^{-}$, too.

For $k=1,2$ and $\mu\in(\mu_{k},\mu_{k+1})$ the equation $l(\varphi)=\mu\varphi$
has $k$ linearly independent bounded solutions $\varphi_{j}(x,\mu)$
$(x\in\mathbb{R}$, $1\leqslant j\leqslant k)$, and the following
asymptotic relations hold \cite{PetKh04}: 
\begin{equation}
\varphi_{j}(x,\mu)=\frac{1}{\sqrt{2\pi}}\sum_{\nu=1}^{1+r_{k}^{+}}\sqrt{|\lambda_{\nu}^{\prime+}(\mu)|}A_{j\nu}^{+}(\mu)e^{ix\lambda_{\nu}^{+}(\mu)}[1+o(1)]\quad\left(x\to\infty\right),
\end{equation}
\[
\varphi_{j}(x,\mu)=\frac{1}{\sqrt{2\pi}}\sum_{\nu=2-r_{k}^{-}}^{2}\sqrt{|\lambda_{\nu}^{\prime-}(\mu)|}A_{j\nu}^{-}(\mu)e^{ix\lambda_{\nu}^{-}(\mu)}[1+o(1)]\quad\left(x\to-\infty\right).
\]
Assuming the notations 
\[
B_{j\nu}(\mu)=\begin{cases}
A_{j\nu}^{+}(\mu) & {\hbox{for}}\quad1\leqslant\nu\leqslant r_{k}^{+}\\
A_{j,\nu-r_{k}^{+}+r_{k}^{-}}^{-}(\mu) & {\hbox{for}}\quad r_{k}^{+}<\nu\leqslant k,
\end{cases}
\]
\[
C_{j\nu}(\mu)=\begin{cases}
A_{j\nu}^{-}(\mu) & {\hbox{for}}\quad1\leqslant\nu\leqslant r_{k}^{-}\\
A_{j,\nu+r_{k}^{+}-r_{k}^{-}}^{+}(\mu) & {\hbox{for}}\quad r_{k}^{-}<\nu\leqslant k
\end{cases}
\]
the matrices 
\begin{equation}
B(\mu)=(B_{j\nu}(\mu))_{j,\nu=1}^{k},\quad C(\mu)=(C_{j\nu}(\mu))_{j,\nu=1}^{k}
\end{equation}
are nondegenerate and satisfy the relation 
\begin{equation}
B(\mu)B^{*}(\mu)=C(\mu)C^{*}(\mu).
\end{equation}
An arbitrary nondegenerate matrix can be taken as one of these matrices,
and this will uniquely determine the other matrix as well as the solutions
$\varphi_{j}(x,\mu)$. From (7), it follows that if one of the matrices
(6) is unitary, then so is the other one, too.
\begin{lem}
For the matrices (6) the following statements are true:
\begin{enumerate}
\item If one of the matrices (6) is measurable (continuous) on the intervals
$(\mu_{1},\mu_{2})$ and $(\mu_{2},\mu_{3})$, then the other is measurable
(continuous) on the same intervals.
\item If the potential $q$ satisfies the condition (4) and one of the matrices
(6) is differentiable (continuously differentiable) on the intervals
$(\mu_{1},\mu_{2})$ and $(\mu_{2},\mu_{3})$, then the other one
is also differentiable (continuously differentiable) on the same intervals.
\end{enumerate}
\end{lem}
\noindent \textit{Proof:} is based on Theorem 1 and explicit formulas
representing the entries of one of the matrices (6) by means of the
entries of the other.
\begin{rem}
\noindent For $r_{k}^{+}=1$, the coefficients $A_{j\nu}^{+}(\mu)$
in (5) are entries of one of the matrices (6). This is not true for
$r_{k}^{+}=0$. However, one can prove that in the latter case measurability,
continuity and, under (4), differentiability and continuous differentiability
of one of the matrices (6) extends to these coefficients, too. The
same is true for the coefficients $A_{j\nu}^{-}(\mu)$.
\end{rem}
Henceforth we assume that $B(\mu)$ and $C(\mu)$ are unitary matrices
and their elements are measurable functions (in particular, one of
them could be the identity matrix). Normalizing in this way, we call
the system of solutions $\varphi_{j}(x,\mu)$ $(\mu\in(\mu_{k},\mu_{k+1})$,
$1\leqslant j\leqslant k)$ the normalized system of generalized eigenfunctions
of the operator $L$, corresponding to the value $\mu$.

Consider the matrices 
\[
A^{+}(\mu)=(A_{j\nu}^{+}(\mu))_{1\leqslant j\leqslant k,\ 1\leqslant\nu\leqslant1+r_{k}^{+}},\quad A^{-}(\mu)=(A_{j\nu}^{-}(\mu))_{1\leqslant j\leqslant k,\ 2-r_{k}^{-}\leqslant\nu\leqslant2}.
\]
For any $x\in\mathbb{R}$ and any $\mu\in(\mu_{k},\mu_{k+1})$ $(k=1,2)$,
$\varphi(x,\mu)$ will denote the vector-column, consisting of the
solutions $\varphi_{j}(x,\mu)$ $(1\leqslant j\leqslant k)$.
\begin{lem}
Let $\varphi_{j}(x,\mu)$ $(\mu\in(\mu_{k},\mu_{k+1}),$ $1\leqslant j\leqslant k)$
be a normalized system of generalized eigenfunctions of the operator
$L$. Then for any measurable unitary matrix $U(\mu)=(U_{ij}(\mu))_{1}^{k}$,
the functions $\widetilde{\varphi}_{j}(x,\mu)$ $(\mu\in(\mu_{k},\mu_{k+1}),$
$1\leqslant j\leqslant k)$ determined by the relation 
\begin{equation}
\widetilde{\varphi}(x,\mu)=U(\mu)\varphi(x,\mu)\quad(\widetilde{\varphi}(x,\mu)=(\widetilde{\varphi}_{j}(x,\mu))_{1\leqslant j\leqslant k}),
\end{equation}
form a normalized system of generalized eigenfunctions of $L$.

Conversely, for any normalized systems $\varphi_{j}(x,\mu)$ and $\widetilde{\varphi}_{j}(x,\mu)$
$(\mu\in(\mu_{k},\mu_{k+1}),$ $1\leqslant j\leqslant k)$ of generalized
eigenfunctions of $L$, there exists a unique measurable, unitary
matrix $U(\mu)=(U_{ij}(\mu))_{i,j=1}^{k}$ satisfying (8). Moreover,
the following equalities hold: 
\[
\widetilde{A}^{\pm}(\mu)=U(\mu)A^{\pm}(\mu),\quad\widetilde{B}(\mu)=U(\mu)B(\mu),\quad\widetilde{C}(\mu)=U(\mu)C(\mu).
\]

\end{lem}
\noindent \textit{Proof:} is simple and directly follows from the
corresponding definitions. \qed

If the operator $L$ has eigenvalues, we consider an orthonormal system
of eigenfunctions\linebreak{}
$\{\psi_{j}$, $j=1,2,\ldots\}$ of $L$, complete in the closure
of the span of all eigenfunctions of $L$. If the point spectrum $T$
of the operator $L$ is finite (countable), then the system of its
eigenfunctions $\psi_{j}$ is finite (countable). Since the eigenvalues
of $L$ are simple, to each eigenvalue of $L$ corresponds exactly
one eigenfunction from $\{\psi_{j}:j=1,2,\ldots\}$.

According to the results of \cite{PetKh04}, any function $f\in L^{2}(\mathbb{R})$
has Fourier expansion 
\begin{equation}
f(x)=\sum_{j}\psi_{j}(x)\int\limits _{-\infty}^{\infty}f(t)\overline{\psi_{j}(t)}\;dt+\sum_{k=1}^{2}\sum_{j=1}^{k}\int\limits _{\mu_{k}}^{\mu_{k+1}}F_{j}(\mu)\varphi_{j}(x,\mu)d\mu\quad(x\in\mathbb{R}),
\end{equation}
and for any $f,g\in L^{2}(\mathbb{R})$ the generalized Parseval equality
is true: 
\begin{equation}
\int\limits _{-\infty}^{\infty}f(x)\overline{g(x)}\;dx=\sum_{j}\int\limits _{-\infty}^{\infty}f(t)\overline{\psi_{j}(t)}\;dt\int\limits _{-\infty}^{\infty}\overline{g(t)}\psi_{j}(t)\;dt+\sum_{k=1}^{2}\sum_{j=1}^{k}\int\limits _{\mu_{k}}^{\mu_{k+1}}F_{j}(\mu)\overline{G_{j}(\mu)}\;d\mu,
\end{equation}
where 
\begin{equation}
F_{j}(\mu)=\int\limits _{-\infty}^{\infty}f(t)\overline{\varphi_{j}(t,\mu)}\;dt,\quad G_{j}(\mu)=\int\limits _{-\infty}^{\infty}g(t)\overline{\varphi_{j}(t,\mu)}\;dt,
\end{equation}
and the last integral in (9) converges in the norm of $L^{2}(\mathbb{R})$,
and for $\mu\in(\mu_{k},\mu_{k+1})$ the integrals in (11) converge
in the norm of the space $L^{2}(\mu_{k},\mu_{k+1})$ (if $L$ has
no eigenvalues, then the first sums in (9) and (10) vanish).

Now we introduce the scattering data for the operator $L$. To this
end, for $\mu\in(\mu_{k},\mu_{k+1})$ and $1\leqslant j,\nu\leqslant1+r_{k}^{+}\;\:(k=1,2)$
we set 
\begin{equation}
S_{j\nu}^{+}(\mu)=\sum_{l=1}^{k}\sqrt{|\lambda_{\nu}^{\prime+}(\mu)|}\sqrt{|\lambda_{j}^{\prime+}(\mu)|}A_{l\nu}^{+}(\mu)\overline{A_{lj}^{+}(\mu)}
\end{equation}
and consider the square matrix $S^{+}(\mu)=(S_{j\,\nu}^{+}(\mu))_{j,\nu=1}^{1+r^{+}(\mu)}$
of order $1+r^{+}(\mu)$. By Lemma 2, this matrix is independent of
the choice of normalized system of generalized eigenfunctions. Hence
Lemma 1 and Remark 1 imply that under condition (1) the matrix-function
$S^{+}$ is continuous and under condition (4) it is continuously
differentiable in the intervals $(\mu_{1},\mu_{2})$ and $(\mu_{2},\mu_{3})$.

Let $T$ be the point spectrum of $L$, $\mu\in T$ and $\psi\left(\cdot,\mu\right)$
be a corresponding normalized eigenfunction. Then the following asymptotic
relation holds: 
\begin{equation}
\psi(x,\mu)=c^{+}(\mu)e^{ix\lambda_{1}^{+}(\mu)}[1+o(1)]\quad{\hbox{as}}\quad x\to\infty,
\end{equation}
where $c^{+}(\mu)$ is a nonzero complex number. One can easily see
that the numbers 
\begin{equation}
N^{+}(\mu)=|c^{+}(\mu)|^{2}\quad(\mu\in T)
\end{equation}
do not depend on the choice of the normalized eigenfunctions $\psi$.

Now consider the data 
\begin{equation}
\{T,\,N^{+}(\mu)\;\:(\mu\in T),\,S^{+}(\mu)\;\:(\mu\in(\mu_{1},\mu_{2})\cup(\mu_{2},\mu_{3}))\}
\end{equation}
called the \emph{right scattering data} of the operator $L$.

\emph{The inverse scattering problem for the operator $L$ consists
in recovery of $L$ from the knowledge of the right scattering data
(15).}

It is known (see \cite{Mar77}, pp. 162--166) that if for each $a\in\mathbb{R}$
the function $q$ satisfies the condition 
\[
\int\limits _{a}^{\infty}(1+|x|)|q(x)-a^{+}|\;dx<\infty
\]
with some number $a^{+}\in\mathbb{R}$, then for any $\lambda\in\mathbb{C},\,{\hbox{\rm Im}}\lambda\geqslant0$
the equation $l(y)=(\lambda^{2}+a^{+})y$ has a solution $y^{+}(x,\lambda)$
for which the representations
\begin{equation}
\begin{gathered}y^{+}(x,\lambda)=e^{i\lambda x}+\int\limits _{x}^{\infty}e^{i\lambda t}K^{+}(x,t)\;dt\quad(-\infty<x<\infty),\hfill\\
e^{i\lambda x}=y^{+}(x,\lambda)+\int\limits _{x}^{\infty}y^{+}(t,\lambda)H^{+}(x,t)\;dt\quad(-\infty<x<\infty)\hfill
\end{gathered}
\end{equation}
hold. Here the kernels $K^{+}(x,t)$ and $H^{+}(x,t)$ $(-\infty<x\leqslant t<\infty)$
do not depend on $\lambda$, are real-valued continuous functions
of two variables $x,\,t$ and satisfy the relations 
\begin{equation}
H^{+}(x,\,t)+K^{+}(x,\,t)+\int\limits _{x}^{t}H^{+}(x,\xi)K^{+}(\xi,\,t)\;d\xi=0,
\end{equation}
\begin{equation}
K^{+}(x,\,t)+H^{+}(x,\,t)+\int\limits _{x}^{t}K^{+}(x,\xi)H^{+}(\xi,\,t)\;d\xi=0.
\end{equation}
The estimates 
\begin{equation}
|K^{+}(x,\,t)|\leqslant\frac{1}{2}h^{+}\left(\frac{x+t}{2}\right)\exp\left[h_{1}^{+}(x)-h_{1}^{+}\left(\frac{x+t}{2}\right)\right],
\end{equation}
\begin{equation}
|H^{+}(x,\,t)|\leqslant\frac{1}{2}h^{+}\left(\frac{x+t}{2}\right)\exp\left\{ h_{1}^{+}(x)-h_{1}^{+}\left(\frac{x+t}{2}\right)+\exp\left[h_{1}^{+}(x)-h_{1}^{+}\left(\frac{x+t}{2}\right)\right]-1\right\} 
\end{equation}
hold, where $-\infty<x\leqslant t<\infty$ and 
\[
h^{+}(x)=\int\limits _{x}^{\infty}|q(t)-a^{+}|\;dt,\quad h_{1}^{+}(x)=\int\limits _{x}^{\infty}h^{+}(t)\;dt.
\]
For $a\in\mathbb{R}$ and $1\leqslant p\leqslant\infty$ the operators
$K_{a}^{+},H_{a}^{+}$, defined by 
\[
(K_{a}^{+}f)(x)=\int\limits _{x}^{\infty}K^{+}(x,\,t)f(t)\;dt\quad(f\in L^{p}(a,\infty),\ x>a),
\]
\[
(H_{a}^{+}f)(x)=\int\limits _{x}^{\infty}H^{+}(x,\,t)f(t)\;dt\quad(f\in L^{p}(a,\infty),\ x>a),
\]
are bounded in $L^{p}(a,\infty)$, $I+K_{a}^{+}$ is invertible and
\begin{equation}
(I+K_{a}^{+})^{-1}=I+H_{a}^{+},
\end{equation}
\begin{equation}
K^{+}(x,x)=\frac{1}{2}\int\limits _{x}^{\infty}(q(t)-a^{+})\;dt\quad\left(x\in\mathbb{R}\right).
\end{equation}

From now on, the condition (4) is assumed to be satisfied.

In analogy with \cite{Khach83}, \cite{Bab89}, we consider the function
\[
\widetilde{F}^{+}(x,t)=\sum_{\mu\in T}\frac{N^{+}(\mu)}{|\lambda_{1}^{+}(\mu)|^{2}}\left(e^{ix\lambda_{1}^{+}(\mu)}-1\right)\left(e^{it\lambda_{1}^{+}(\mu)}-1\right)+
\]
\begin{equation}
+\frac{1}{2\pi}\int\limits _{\mu_{1}}^{\infty}\,\sum_{\nu,j=1}^{1+r^{+}(\mu)}\frac{S_{j\nu}^{+}(\mu)}{\lambda_{\nu}^{+}(\mu)\overline{\lambda_{j}^{+}(\mu)}}\bigg(e^{ix\lambda_{\nu}^{+}(\mu)}-1\bigg)\bigg(e^{\overline{it\lambda_{j}^{+}(\mu)}}-1\bigg)\;d\mu-\omega(x,t)\quad(x,\,t\in\mathbb{R}),
\end{equation}
where 
\[
\omega(x,t)=\begin{cases}
\min\{|x|,\,|t|\} & \text{ for}\quad xt\geqslant0.\\
0 & \text{ for}\quad xt<0.
\end{cases}
\]
According to what was said above, the point spectrum $T$ of the operator
$L$ is a finite set. Therefore the first sum in (23) contains a finite
number of summands. Since $T\subset(-\infty,\mu_{1})$, for $\mu\in T$
the numbers $\lambda_{1}^{+}(\mu)$ lie in the upper part of the imaginary
axis. Hence the mentioned sum is a real number. We will see below
that the integral in (23) is convergent (in the usual sense).
\begin{thm}
The derivative 
\[
F^{+}(x,t)=\frac{\partial^{2}\widetilde{F}^{+}(x,t)}{\partial x\partial t}\quad\left(x,t\in\mathbb{R}\right)
\]
exists and is continuous, real-valued and symmetric: 
\begin{equation}
F^{+}(x,t)=F^{+}(t,x)\quad(x,t\in\mathbb{R}).
\end{equation}
Moreover, the following equalities hold: 
\begin{equation}
F^{+}(x,t)=H^{+}(x,t)+\int\limits _{t}^{\infty}H^{+}(x,\xi)H^{+}(t,\xi)\;d\xi\quad(-\infty<x\leqslant t<\infty),
\end{equation}
\begin{equation}
F^{+}(x,t)=H^{+}(t,x)+\int\limits _{x}^{\infty}H^{+}(x,\xi)H^{+}(t,\xi)\;d\xi\quad(-\infty<t\leqslant x<\infty).
\end{equation}

\end{thm}
\noindent \textit{Proof:} For $\mu\in(\mu_{k},\mu_{k+1})$, $1\leqslant l\leqslant k$
$(k=1,2)$ we denote 
\begin{equation}
F_{N}^{+}(x,t,\mu)=N^{+}(\mu)e^{ix\lambda_{1}^{+}(\mu)}e^{it\lambda_{1}^{+}(\mu)}\quad(\mu\in T),
\end{equation}
\begin{equation}
v_{l}(x,\mu)=\frac{1}{\sqrt{2\pi}}\sum_{\nu=1}^{1+r^{+}(\mu)}\sqrt{|\lambda_{\nu}^{'+}(\mu)|}A_{l\nu}^{+}(\mu)e^{ix\lambda_{\nu}^{+}(\mu)}.
\end{equation}
The following equalities hold 
\[
\frac{N^{+}(\mu)}{|\lambda_{1}^{+}(\mu)|^{2}}\left(e^{ix\lambda_{1}^{+}(\mu)}-1\right)\left(e^{it\lambda_{1}^{+}(\mu)}-1\right)=\int\limits _{0}^{x}\!\!\int\limits _{0}^{t}F_{N}^{+}(\xi,\eta,\mu)\;d\eta\;d\xi\quad(\mu\in T),
\]
\[
\frac{1}{2\pi}\sum_{\nu,j=1}^{1+r^{+}(\mu)}\frac{S_{j\nu}^{+}(\mu)}{\lambda_{\nu}^{+}(\mu)\overline{\lambda_{j}^{+}(\mu)}}\bigg(e^{ix\lambda_{\nu}^{+}(\mu)}-1\bigg)\bigg(e^{\overline{it\lambda_{j}^{+}(\mu)}}-1\bigg)=
\]
\[
=\sum_{l=1}^{k}\int\limits _{0}^{x}v_{l}(\xi,\mu)\;d\xi\int\limits _{0}^{t}\overline{v_{l}(\eta,\mu)}\;d\eta\quad(\mu\in(\mu_{k},\mu_{k+1}),\ k=1,2)
\]
(to obtain the first equality, it is necessary to take into account
that $|\lambda_{1}^{+}(\mu)|^{2}=-[\lambda_{1}^{+}(\mu)]^{2}$ since
$i\lambda_{1}^{+}(\mu)$ are real numbers). Consequently, formula
(23) takes the form 
\begin{multline}
\widetilde{F}^{+}(x,t)=\sum_{\mu\in T}\int\limits _{0}^{x}\int\limits _{0}^{t}F_{N}^{+}(\xi,\eta,\mu)\;d\eta\;d\xi\\
+\sum_{k=1}^{2}\sum_{l=1}^{k}\int\limits _{\mu_{k}}^{\mu_{k+1}}\left[\int\limits _{0}^{x}v_{l}(\xi,\mu)\;d\xi\int\limits _{0}^{t}\overline{v_{l}(\eta,\mu)}\;d\eta\right]d\mu-\omega(x,t).\qquad
\end{multline}
(13) implies the equality 
\begin{equation}
\psi(x,\mu)=c^{+}(\mu)y^{+}(x,\lambda_{1}^{+}(\mu)).
\end{equation}
Since the function $y^{+}(x,\lambda_{1}^{+}(\mu))$ is real-valued,
we have 
\begin{equation}
\overline{\psi(x,\mu)}=\overline{c^{+}(\mu)}y^{+}(x,\lambda_{1}^{+}(\mu)).
\end{equation}
Further, in view of (27) and (16), we have 
\[
F_{N}^{+}(x,t,\mu)=N^{+}(\mu)\left[y^{+}(x,\lambda_{1}^{+}(\mu))+\int\limits _{x}^{\infty}y^{+}(\xi,\lambda_{1}^{+}(\mu))H^{+}(x,\xi)\;d\xi\right]\times
\]
\[
\times\left[y^{+}(t,\lambda_{1}^{+}(\mu))+\int\limits _{t}^{\infty}y^{+}(\eta,\lambda_{1}^{+}(\mu))H^{+}(t,\eta)\;d\eta\right].
\]
From this equality and (14), (30), (31) we obtain 
\begin{equation}
F_{N}^{+}(x,t,\mu)=\left[\psi(x,\mu)+\int\limits _{x}^{\infty}\psi(\xi,\mu)H^{+}(x,\xi)\;d\xi\right]\left[\overline{\psi(t,\mu)}+\int\limits _{t}^{\infty}\overline{\psi(\eta,\mu)}H^{+}(t,\eta)\;d\eta\right].
\end{equation}
(5) implies that 
\[
\varphi_{l}(x,\mu)=\frac{1}{\sqrt{2\pi}}\sum_{\nu=1}^{1+r^{+}(\mu)}\sqrt{|{\lambda'_{\nu}}^{+}(\mu)|}A_{l\nu}^{+}(\mu)y^{+}(x,\lambda_{\nu}^{+}(\mu)),
\]
hence, using (28) and (16), we get 
\begin{equation}
v_{l}(x,\mu)=\varphi_{l}(x,\mu)+\int\limits _{x}^{\infty}\varphi_{l}(\xi,\mu)H^{+}(x,\xi)\;d\xi\quad(1\leqslant l\leqslant k,\quad\mu\in(\mu_{k},\mu_{k+1})).
\end{equation}
Consider the function $f_{x}\left(\cdot\right)\;\:\left(x\in\mathbb{R}\right)$
defined by the formula 
\[
f_{x}(\xi)=\begin{cases}
0 & {\hbox{if}}\quad\xi<0\\
{\displaystyle {1+\int\limits _{0}^{\xi}H^{+}(\eta,\xi)\;d\eta}} & {\hbox{if}}\quad0\leqslant\xi<x\\
{\displaystyle {\int\limits _{0}^{x}H^{+}(\eta,\xi)\;d\eta}} & {\hbox{if}}\quad\xi\geqslant x
\end{cases}
\]
for $x\geqslant0$ and by the formula\pagebreak{}

\noindent 
\[
f_{x}(\xi)=\begin{cases}
0 & {\hbox{if}}\quad\xi<x\\
{\displaystyle {-1-\int\limits _{x}^{\xi}H^{+}(\eta,\xi)\;d\eta}} & {\hbox{if}}\quad x\leqslant\xi<0\\
{\displaystyle {\int\limits _{0}^{x}H^{+}(\eta,\xi)d\eta}} & {\hbox{if}}\quad\xi\geqslant0
\end{cases}
\]
for $x<0$. The integration of (32) and (33) yields 
\[
\int\limits _{0}^{x}\!\int\limits _{0}^{t}F_{N}^{+}(\xi,\eta,\mu)\;d\eta\;d\xi=\int\limits _{-\infty}^{\infty}\psi(\xi,\mu)f_{x}(\xi)\;d\xi\int\limits _{-\infty}^{\infty}\overline{\psi(\xi,\mu)}f_{t}(\xi)\;d\xi,
\]
\[
\int\limits _{0}^{x}v_{l}(\xi,\mu)\;d\xi=\int\limits _{-\infty}^{\infty}\varphi_{l}(\xi,\mu)f_{x}(\xi)\;d\xi.
\]
In view of these equalities, (29) takes the form 
\[
\widetilde{F}^{+}(x,t)=\sum_{\mu\in T}\int\limits _{-\infty}^{\infty}\psi(\xi,\mu)f_{x}(\xi)\;d\xi\int\limits _{-\infty}^{\infty}\overline{\psi(\xi,\mu)}f_{t}(\xi)\;d\xi+
\]
\begin{equation}
+\sum_{k=1}^{2}\sum_{l=1}^{k}\int\limits _{\mu_{k}}^{\mu_{k+1}}\left[\int\limits _{-\infty}^{\infty}\varphi_{l}(\xi,\mu)f_{x}(\xi)\;d\xi\int\limits _{-\infty}^{\infty}\overline{\varphi_{l}(\xi,\mu)}f_{t}(\xi)\;d\xi\right]\;d\mu-\omega(x,t).
\end{equation}
By the estimate (20), $f_{x}(.)\in L^{2}(\mathbb{R})$ for any $x\in\mathbb{R}$.
Therefore, by (34) and the generalized Parseval equality (10), we
get 
\[
\widetilde{F}^{+}(x,t)=\int\limits _{-\infty}^{\infty}f_{x}(\xi)f_{t}(\xi)\;d\xi-\omega(x,t).
\]
This equality implies the existence of the continuous mixed derivative
$F^{+}(x,t)=\frac{\partial^{2}\widetilde{F}^{+}(x,t)}{\partial x\partial t}$
on $\mathbb{R}^{2}$ and the representations (25), (26). These representations
show that $F^{+}$ is real-valued and symmetric. The proof is complete.
\begin{thm}
\noindent The function $F^{+}(x,t)$ and the kernel $K^{+}(x,t)$
satisfy 
\begin{equation}
F^{+}(x,t)+K^{+}(x,t)+\int\limits _{x}^{\infty}K^{+}(x,\xi)F^{+}(\xi,t)\;d\xi=0\quad(-\infty<x\leqslant t<\infty),
\end{equation}
\begin{equation}
F^{+}(x,t)+\int\limits _{x}^{\infty}K^{+}(x,\xi)F^{+}(\xi,t)\;d\xi=H^{+}(t,x)\quad(-\infty<t\leqslant x<\infty).
\end{equation}

\end{thm}
\noindent \textit{Proof:} For a fixed $a\in\mathbb{R}$, in $L^{2}(a,\infty)$
consider the operators $K_{a}^{+}$, $H_{a}^{+}$ and the operator
$F_{a}^{+}$ defined by 
\[
(F_{a}^{+}f)(x)=\int\limits _{a}^{\infty}F^{+}(x,t)f(t)\;dt\quad(f\in L^{2}(a,\infty),\ x>a).
\]
(25), (26) and (20) show that there exists a decreasing summable function
$\sigma$ on $[a,\infty)$ such that 
\begin{equation}
|F^{+}(x,t)|\leqslant\sigma\left(\frac{x+t}{2}\right)\quad(a\leqslant x,t<\infty).
\end{equation}
Hence the operator $F_{a}^{+}$ is bounded. Moreover, by (25) and
(26), 
\begin{equation}
I+F_{a}^{+}=(I+H_{a}^{+})(I+H_{a}^{*+}),
\end{equation}
and by (38) and (21), 
\[
K_{a}^{+}+F_{a}^{+}+K_{a}^{+}F_{a}^{+}=H_{a}^{*+}.
\]
For the corresponding kernels we get 
\begin{equation}
F^{+}(x,t)+K^{+}(x,t)+\int\limits _{x}^{\infty}K^{+}(x,\xi)F^{+}(\xi,t)\;d\xi=0\quad(a<x\leqslant t<\infty),
\end{equation}
\begin{equation}
F^{+}(x,t)+\int\limits _{x}^{\infty}K^{+}(x,\xi)F^{+}(\xi,t)\;d\xi=H^{+}(t,x)\quad(a<t\leqslant x<\infty).
\end{equation}
Since $a$ was arbitrary, the equalities (35) and (36) follow from
(39) and (40). The proof is complete.
\begin{lem}
\noindent For any $p\in\left[1,\infty\right]$ and $x\in\mathbb{R}$
the function $K^{+}(x,\cdot)$ is the unique solution of the equation
(35) in the space $L^{p}(x,\infty)$.
\end{lem}
\noindent \textit{Proof:} Fix $p\in\left[1,\infty\right]$ and $x\in\mathbb{R}$.
Consider the following integral operators on $L^{p}(0,\infty)$: 
\[
(K_{x}f)(\xi)=\int\limits _{\xi}^{\infty}K^{+}(x+\xi,x+\eta)f(\eta)\;d\eta\quad(\xi>0),
\]
\[
(K_{x}^{*}f)(\xi)=\int\limits _{0}^{\xi}K^{+}(x+\eta,x+\xi)f(\eta)\;d\eta\quad(\xi>0),
\]
\[
(H_{x}f)(\xi)=\int\limits _{\xi}^{\infty}H^{+}(x+\xi,x+\eta)f(\eta)\;d\eta\quad(\xi>0),
\]
\[
(H_{x}^{*}f)(\xi)=\int\limits _{0}^{\xi}H^{+}(x+\eta,x+\xi)f(\eta)\;d\eta\quad(\xi>0),
\]
\[
(G_{x}f)(\xi)=\int\limits _{0}^{\infty}F^{+}(x+\xi,x+\eta)f(\eta)\;d\eta\quad(\xi>0).
\]
By the estimates (19), (20) and (37) for $K^{+}(x,t)$, $H^{+}(x,t)$
and $F^{+}(x,t)$, the operators $K_{x}$, $K_{x}^{*}$, $H_{x}$,
$H_{x}^{*}$ and $G_{x}$ are bounded. It is obvious that $K_{x}^{*}$
and $H_{x}^{*}$ are the conjugate operators of $K_{x}$ and $H_{x}$
in the case $p=2$.

The equalities (17) and (18) imply 
\[
(I+H_{x})(I+K_{x})=I,\quad(I+K_{x})(I+H_{x})=I,
\]
i.e., the operator $I+H_{x}$ is invertible and 
\begin{equation}
(I+H_{x})^{-1}=I+K_{x}.
\end{equation}
Similarly, $I+H_{x}^{*}$ is an invertible operator and 
\begin{equation}
(I+H_{x}^{*})^{-1}=I+K_{x}^{*}.
\end{equation}
It follows from (25) and (26) that 
\[
I+G_{x}=(I+H_{x})(I+H_{x}^{*}).
\]
Using this and (41), (42), we conclude that $I+G_{x}$ is an invertible
operator.

To show that equation (35) can have no more that one solution in $L^{p}(x,\infty)$,
suppose that 
\[
g_{1}(t)+F^{+}(x,t)+\int\limits _{x}^{\infty}g_{1}(\xi)F^{+}(\xi,t)\;d\xi=0\quad(x<t<\infty),
\]
\[
g_{2}(t)+F^{+}(x,t)+\int\limits _{x}^{\infty}g_{2}(\xi)F^{+}(\xi,t)\;d\xi=0\quad(x<t<\infty),
\]
for some functions $g_{j}\in L^{p}(x,\infty)$ $(j=1,2)$. Subtract
the second equation from the first one to obtain 
\begin{equation}
g_{1}(t)-g_{2}(t)+\int\limits _{x}^{\infty}[g_{1}(\xi)-g_{2}(\xi)]F^{+}(\xi,t)\;d\xi=0\quad(x<t<\infty).
\end{equation}
Put $g(u)=g_{1}(u+x)-g_{2}(u+x)$ $(u>0)$. Using (24), we bring (43)
to the form 
\[
g(t-x)+\int\limits _{x}^{\infty}F^{+}(t,\xi)g(\xi-x)d\xi=0\quad(x<t<\infty).
\]
Substituting $u=t-x$ and $v=\xi-x$, we get 
\[
g(u)+\int\limits _{0}^{\infty}F^{+}(u+x,v+x)g(v)dv=0\quad(0<u<\infty)
\]
which is the same as 
\[
(I+G_{x})g=0.
\]
By the invertability of the operator $I+G_{x}$, this means that $g(u)=0$
$(0<u<\infty)$, i.e., 
\[
g_{1}(u+x)-g_{2}(u+x)=0\quad(0<u<\infty).
\]
Taking $u=t-x$ $(x<t<\infty)$, we complete the proof: 
\[
g_{1}(t)=g_{2}(t)\quad(x<t<\infty).
\]

The above lemma shows that the kernel $K^{+}(x,t)$ can be uniquely
recovered by the right scattering data. This implies that the potential
$q$ can also be uniquely recovered by the right scattering data.
Indeed, by the unitarity of the matrices (6) and (12), for $\mu>\mu_{1}$
\[
S_{11}^{+}(\mu)=\lambda_{1}^{\prime+}(\mu)=\frac{1}{2\sqrt{\mu-a^{+}}}.
\]
Therefore the constant $a^{+}$ is determined by the right scattering
data: 
\[
a^{+}=\lim_{\mu\to\infty}\left\{ \mu-\frac{1}{4[S_{11}^{+}(\mu)]^{2}}\right\} .
\]
Moreover, (22) shows that 
\[
q(x)=a^{+}-2\frac{d}{dx}K^{+}(x,x),
\]
and hence the potential $q$ is uniquely recovered by the right scattering
data. Practically, the recovery of $q$ reduces to solution of the
linear integral equation (35).

The author express deep gratitude to Prof. I. G. Khachatryan for posing
the problem and useful discussions of the obtained results.

\end{document}